\documentclass[12pt]{article}

\usepackage{amssymb,amsmath,latexsym}
\usepackage[dvips]{graphics}
\usepackage{mathrsfs}
\usepackage{color}
\usepackage{comment}

\newtheorem{thm}{Theorem}[section]
\newtheorem{prop}[thm]{Proposition}

\newtheorem{lemma}[thm]{Lemma}
\newtheorem{remark}[thm]{Remark}

\newtheorem{cor}[thm]{Corollary}
\numberwithin{equation}{section}

\definecolor{dgreen}{rgb}{0, 0.6, 0.1}

 \def\sB {{\mathcal B}}

 \def\bN {{\mathbb N}} 
\def\bP {{\mathbb P}}  \def\bR {{\mathbb R}}
 \def\bT {{\mathbb T}} 
  
 \def\bZ {{\mathbb Z}}

\def\ZZ{\mathbb{Z}}

\def\<{\langle}
\def\>{\rangle}
\def\pf{\noindent{\bf Proof.} }
\def\proof{\noindent{\bf Proof.} }

\def\E{{\bf E}}

\def\P{{\bP}}

\def\qed{{\hfill $\Box$\medskip}}

\def\to{\rightarrow}

\def\epsilon{\varepsilon}

\allowdisplaybreaks[4]

\def\be{\begin{equation}}

\def\ee{\end{equation}}

\def\al{\alpha}

\begin{document}
\title{\bf Gaussian bounds and  Collisions of variable speed random walks on  lattices with power law  conductances}

\author{
 Xinxing Chen}

\maketitle

\begin{abstract}
We consider a  weighted lattice $\bZ^d$ with conductance
$\mu_e=|e|^{-\al}$.   We show that the heat kernel of a variable
speed  random walk on it satisfies a two-sided Gaussian bound by
using an intrinsic metric. We also show that when $d=2$ and $\al\in
(-1,0)$, two independent random walks on such weighted lattice
 will collide infinite many times while they are transient.

\end{abstract} \noindent{\bf 2000 MR subject
classification:}   60G50, 58J35


\section{Introduction} \label{s:intro} 
In \cite{HS}, Hebisch and Saloff-Coste proved  that  when  a group
has polynomial volume growth of order $D$, the heat kernel of a
constant speed random walk on the group satisfies a two-sided
Gaussian estimate, i.e.,
$$
c_1t^{-D/2}\exp\left(-c_2\frac{\rho(x,y)^2}{t}\right) \le
p_t(x,y)\le c_3t^{-D/2}\exp\left(-c_4\frac{\rho(x,y)^2}{t}\right).
$$
where $\rho(x,y)$ is a metric on the group.
  Delmotte \cite{DT}  gave    equivalence of Gaussian
bounds, parabolic Harnack inequalities, and the combination of
volume regularity and Poincar\'e inequality.  Later, there are many
papers, such as  \cite{ABDH,MB1, BBK,BC, BD, Sap}, showing that
Gaussian bounds hold for lattice  $\bZ^d$ with different  random
conductances.
%
In this paper, we consider a deterministic weighted lattice which
does not satisfy  Poincar$\acute{e}$ inequalities for all
(sufficiently large) balls or volume doubling property, show that a
variable random walk on it also satisfies the two-side Gaussian
bound, but with a
 metric which is not comparable to the Euclidean metric.


 Let $\al\in \bR$. For
$x,y\in \bZ^d$ with $|x-y|_1=1$, we set $\mu_{xy}=(| x|_\infty \vee
| y|_\infty)^{-\al}$ for the  conductance of $(x,y)$. For
convenience, we set $\mu_{xy}=0$ if
 $x$ and $ y$ are not nearest neighbor. Write $\mu_x=\sum_y\mu_{xy}$ and $\nu_x=(|x|_\infty\vee 1)^\al$ for each $x\in \bZ^d$.  Let $X=\{X_t:
t\ge 0\}$ be a continuous time random walk
 on the lattice $\bZ^d$ with generator
$$ \mathscr{L}f(x) = \frac{1}{\nu_x}\sum_{y\in \ZZ^d} (f(y)-f(x) )\mu_{xy}.$$
Then $X$ is a variable speed random walk waiting for an
exponentially
 distributed time with mean $\frac{\nu_x}{\mu_x}\asymp |x|_\infty^{2\al}$ before jumping.
  $~$ The
transition density of $X$  with respect to $\nu$ is denoted by
$$p_t(x,y)=\frac{\P_x(X_t=y)}{\nu_y}.$$
To show the  Gaussian bounds hold,  we introduce a  metric $\rho$ of
$\bZ^d$.
 We call $x_0\cdots x_m$  a path if
$|x_{i+1}-x_i|_1=1$ for each $i<m$. Let  $\rho(x,x)=0$ for $x\in
\bZ^d$, and for $x,y\in \ZZ^d$ with $y\not=x$ set
$$
\rho(x,y)=\min\left\{\sum_{i=0}^{m}\nu_{z_i}:~z_0z_1\cdots
z_m\text{~is a path~with~}z_0=x~{\rm and~}z_m=y \right\}.
$$
Then there exists a constant $C=C(\al, d)$, such that
\begin{equation}\label{e:1013}
\frac{1}{\nu_x}\sum_{y\sim x}\rho(x,y)^2\mu_{xy}\le
C~~{\rm~for~all~}x.
\end{equation}
Metrics satisfying (\ref{e:1013}) are called  intrinsic metrics, see
 \cite{FLW,SK}. One may expect that analogues of diffusion
processes on manifolds hold using the intrinsic metrics for random
walks on graphs.
 For  $x\in \ZZ^d$ and $r\in \bR^+$, write
$B_{\rho}(x,r)=\{y\in \ZZ^d: \rho(x,y)\le r\}$ for a $\rho-$ball. We
extend $\nu$  to a  measure on $\bZ^d$  and set
$$
V_{\rho}(x,r)=\nu(B_{\rho}(x,r)).
$$

\begin{thm}\label{T:lower1848}
 Let $\al>-1$. Let  $x,y\in \bZ^d$ and $t>0$. If $t< (\nu_x\vee \nu_y)\rho(x,y)$, then
\begin{equation}\label{e:1843w}
p_t(x,y)\le c_1(\nu_x\nu_y)^{-1/2}
\exp\left(-\frac{c_2\rho(x,y)}{\nu_x\vee \nu_y}\Big(1\vee \log
(\frac{(\nu_x\vee \nu_y)\rho(x,y)}{t})\Big) \right).
\end{equation}
 If $t\ge (\nu_x\vee
\nu_y)\rho(x,y)$, then
\begin{equation}\label{e:1845w}
p_t(x,y)\le \frac{c_3}{\sqrt{V_{\rho}(x, t^{1/2})V_{\rho}(y,
t^{1/2})} }\exp\left(-c_4\frac{\rho(x,y)^2}{t}\right)
\end{equation}
and
\begin{equation}\label{lower:1608w}
p_t(x,y)\ge  \frac{c_5}{\sqrt{V_{\rho}(x, t^{1/2})V_{\rho}(y,
t^{1/2})}}\exp\left(-c_6\frac{\rho(x,y)^2}{t}\right).
\end{equation}

\end{thm}

\begin{remark} 
(1) In Lemmas  \ref{c:0822} and  \ref{l:0900}, we give the bounds of
$\rho(x,y)$ and $V_{\rho}(x,t^{1/2})$, respectively.\\
(2) Note that if $\al<-1$ then $\sup_{x,y}\rho(x,y)<\infty$ and $X$
will explode in a finite time.  However, we still do not know
whether the heat kernel of $X$ has Gaussian bounds at the critical
point $\al=-1$.
\end{remark}

 Next, we  consider the collision problem of random walks on these weighted lattices. As
usual, we say that two walks $X$ and $X'$ collide infinitely often
if almost surely there exists a sequence of (random) times
$\{t_i:i\ge 1\}$ with $\lim_it_i=\infty$ such that
$X_{t_i}=X'_{t_i}$ for all $i$. In \cite{Polya}, P$\acute{o}$lya
first studied whether two independent simple random walks on $\bZ^d$
collide infinitely often. He reduced it to the problem of a single
walker returning to his starting point. Later Jain and Pruitt in
\cite{JP} showed the Hausdroff dimension of the intersection of two
independent stable processes, and  Shieh in \cite{SNR} gave a
sufficient condition for infinitely collisions of L$\acute{e}$vy
processes in $\bR$. However, if the walks   are not on a homogeneous
space, the problem will be complicated. Recently in \cite{HP},
Hutchcroft and Peres  use the Mass-Transport Principle to prove that
a recurrent reversible random rooted graph has the infinite
collision property.
 Examples that two recurrent random walks will never
meet, were  shown  in \cite{BF,BPS,KP}. Here, we  give another
example that two transient  random walks will collision infinite
often.

\begin{thm}\label{t:main} Let $\al>-1$. Let $X'$ be an independent copy of $X$.\\
(1) Process $X$ is recurrent if and only if  $\al\ge d-2$.\\
(2) If $d\le 2$, then  $X$ and $X'$  collide infinitely often.\\
(3) If $d\ge 3$, then  $X$ and $X'$  collide finitely often.
\end{thm}

\begin{remark} It is   much interesting  that $X$ is not recurrent
while  $X$ and $X'$  collide  infinitely often
 when $d= 2$ and $\al\in(-1,0)$. Similarly, when $d\ge 3$ and
 $\al\ge d-2$,  $X$ is
 recurrent while  $X$ and $X'$  collide  finitely often.
\end{remark}

 In Section 2, we  obtain some geometric properties of
 the weighted lattice $\bZ^d$. In Section 3, we obtain an upper bound on $p_T(w,w)$ by
 using the approach of Barlow and Chen \cite{BC}, which in turn is based on \cite{KZ,MB1}. In Section 4, we  obtain
  the lower bounds of near diagonal transition probability by  using the result of Delmette \cite{DT}
 directly and a chain argument. In Section 5, we give the proof of Theorem
 \ref{T:lower1848}.
 Section 6 deals with the proof of   Theorem \ref{t:main} by  the two-sided Gaussian bounds.

Throughout this paper, we use the notation $c$, $c'$ etc to denote
fixed positive constants which may vary on each appearance, and
$c_i$  to denote positive constants which are fixed in each
argument. If we need to refer to constant $c_1$ of Lemma 2.1
elsewhere we will use the notation $c_{2.1.1}$. For any two
functions $f$ and $g$,  we say $f\asymp g$  if  there exists
$c_i(\al,d)>0$ such that $ c_1 f\le g \le c_2f.$ For brevity, we
write $|\cdot|_p$ for the $L^p-$norm of the Euclidean space
$\mathbb{R}^d$, while  $|\cdot|$ instead of $|\cdot|_\infty$ for the
$L^\infty-$norm. Write $B(x,r)=\{ y\in \ZZ^d: |y-x|\le r\}$ for an
$L^\infty-$ball.

\section{Some geometric properties}
Fix $\al>-1$ henceforth.
  In this section, we shall estimate the metric
$\rho(x,y)$ and the volume $V_\rho(x,r)$, and give
Poincar$\acute{e}$ inequalities. Let us begin with  the volume of a
path.
\begin{lemma}\label{l:vip0926} Let  $z_0 \cdots z_n$ be a path with $\max\{|z_0|,|z_n|, |z_0-z_n|_1\}\ge n\ge 1$.
Then
\begin{equation}\label{e:1741s}
c_1 n(|z_0|\vee |z_n|)^\al\le \sum_{i=0}^n\nu_{z_i}\le c_2
n(|z_0|\vee |z_n|)^\al.
\end{equation}
\end{lemma}
\pf Without loss generality, we may assume that $|z_0|\ge |z_n|$ in
the following. (Otherwise, relabel $z_{n-k}$ with $z_k$ for all
$k$.) Then
$$|z_0-z_n|_1\le d|z_0-z_n|\le d|z_0|+d|z_n|\le 2d|z_0|.$$ Using the
condition $\max\{|z_0|,|z_n|, |z_0-z_n|_1\}\ge n\ge 1$, we get
\begin{equation}\label{e:0717s}
|z_0|\ge \tfrac{n}{2d}\vee 1.
\end{equation}
Since $z_0\cdots z_n$ is a path, we have $|z_i-z_0|\le i$ for each
$i$. So, $\nu_{z_i}=(|z_i|\vee 1)^\al$ takes value between
$(|z_0|+i)^\al$ and $((|z_0|-i)\vee 1)^\al.$ Hence $ \nu_{z_i}\ge
c|z_0|^\al{\rm~for~}i\le \tfrac{n}{4d}, $ which implies
\begin{equation}\label{e:1813s}
\sum_{i=0}^n\nu_{z_i}\ge \sum_{i\le n/(4d)}\nu_{z_i}\ge c\lceil
n/4d\rceil \nu_{z_0}\ge c' n |z_0|^\al=c'n(|z_0|\vee |z_n|)^\al.
\end{equation}
We have proved the lower bound of   (\ref{e:1741s}). For the upper
bound,  we consider two cases.

Case I: $|z_0|\ge  |z_n|\vee n$. Directly calculate
\begin{align*}
\sum_{i=0}^n\nu_{z_i}\le& \sum_{i=0}^n (
(|z_0|+i)^\al+((|z_0|-i)\vee
1)^\al)\le 2\sum_{i=|z_0|-n}^{i=|z_0|+n} (i\vee 1)^\al\\
\le& c_1\int_{|z_0|-n}^{|z_0|+n} x^\al
dx=\frac{c_1}{1+\al}((|z_0|+n)^{\al+1}-(|z_0|-n)^{\al+1}).
\end{align*}
Since  $\lim_{t\rightarrow
0+}((1+t)^{\al+1}-(1-t)^{\al+1})t^{-1}=2(\al+1)$, we obtain
$$
\sup_{t\in(0,1]}|((1+t)^{\al+1}-(1-t)^{\al+1})t^{-1}|\le c_2.
$$
Substituting  $t=\tfrac{n}{|z_0|}\le 1$ into the above inequality
gives
$$
\sum_{i=0}^n\nu_{z_i}\le
\frac{c_1}{1+\al}((|z_0|+n)^{\al+1}-(|z_0|-n)^{\al+1})\le
\frac{c_1c_2}{1+\al}n|z_0|^{\al}=cn(|z_0|\vee |z_n|)^\al.
$$

Case II: $|z_n|\le |z_0|< n$ and $|z_0-z_n|_1=n$.  Then $z_0\cdots
z_n$ is an $L^1-$geodesic, which implies $\{z_0,\cdots,z_n\}\subset
B(0,n)$ and $ |\{i:z_i\in B(0, r)\}|\le 2dr$  for each $r$. Write
$$k=\lceil \log_2 n\rceil,~~T_0=B(0,1) {\rm~~ and~~} T_l=B(0,
2^{l})-B(0,2^{l-1}) {\rm~~ for~~} l\ge 1.$$ Then
\begin{align*}
\sum_{i=0}^n\nu_{z_i}=&\sum_{l=0}^{k} \sum_{i:z_i\in
T_l}\nu_{z_i}\le c\sum_{l=0}^{k} 2^{l\al}|\{i:z_i\in T_l\}|\le
c\sum_{l=0}^{k}
2^{\al l}|\{i:z_i\in B(0, 2^{l})\}|\\
\le&c\sum_{l=0}^{k}2^{\al
l}(2d\cdot2^{l})=2dc\sum_{l=0}^{k}2^{(\al+1) l}\le c' 2^{(1+\al)
(k-1)}\le c' n^{1+\al}.
\end{align*}
Since  $\tfrac{n}{2d}\le |z_0|<n$, we still have
$\sum_{i=0}^n\nu_{z_i}\le c_2 n(|z_0|\vee |z_n|)^\al$ and prove the
lemma.
\qed\\

For $x\in \ZZ^d$ and  $r\in\bR^+$, we  set
\begin{equation}\label{e:2024w}
\rho_x(r)=(|x|\vee r)^\al r.
\end{equation}
Then $\rho_x(\cdot)$ is strictly increasing and
\begin{equation}\label{e:0935}\left(\frac{r}{s}\right)^{c_1}\le \frac{\rho_x(r)}{\rho_x(s)} \le
\left(\frac{r}{s}\right)^{c_2},~~~~\forall ~~r\ge s>0.\end{equation}
 A simple calculation gives, if $x,y\in
\ZZ^d$ and $r\ge \kappa|x-y|$, then there exists $C=C(\al,\kappa)>0$
such that
\begin{equation}\label{e:1650}
C^{-1} \rho_y(r)\le \rho_x(r)\le C \rho_y(r).
\end{equation}
Set $\rho_x^{-1}(r)= (|x|\vee r^{1/(1+\al)})^{-\al} r$, which is
 the inverse function of $\rho_x$. Then $\rho_x^{-1}(\cdot)$ also
 satisfies (\ref{e:0935}) and (\ref{e:1650}).
\begin{lemma}\label{c:0822} Let $x,y\in \ZZ^d$. Let  $\gamma$ be an $L^1-$geodesic path  from $x$ to
$y$.  Then \\
 $$ \left\{\rho(x,y), \sum_{u\in V(\gamma)}\nu_u,  \sum_{(u,v)\in E(\gamma)}\mu_{uv}^{-1}\right\} \subset [c_1\rho_x(|x-y|),  c_1^{-1}\rho_x(|x-y|)].$$
\end{lemma}
\pf   By (\ref{e:1650}), we have $ \rho_x(|x-y|)\asymp
\rho_y(|x-y|)$. So, we may assume $|x|\ge |y|$ without loss
generality. (Otherwise, exchange $y$ with $x$.) Hence $|x|\ge
\tfrac12(|x|+|y|)\ge \tfrac12|x-y|,$ which implies
\begin{equation}\label{e:0911s}\rho_x(|x-y|)\asymp |x-y|\cdot|x|^\al.\end{equation} Let
$z_0z_1 \cdots z_m $ be   a $\rho-$geodesic path with $z_0=x$ and
$z_m=y$, then by Lemma \ref{l:vip0926},
$$
\rho(x,y)\ge \tfrac{1}{2}\sum_{k=0}^{\lceil |x-y|/2\rceil
}\nu_{z_k}\ge c\lceil |x-y|/2\rceil |x|^\al.
$$
By the definition of $\rho(x,y)$, it is clear that $\sum_{u\in
V(\gamma)}\nu_u\ge \rho(x,y)$. Moreover, by Lemma \ref{l:vip0926},
$$
\sum_{u\in V(\gamma)}\nu_u\le c|x-y|_1(|x|\vee|y|)^\al\le 2d
c|x-y||x|^\al.
$$
Since $\mu_{uv}^{-1}\asymp (|u|\vee 1)^{\al}= \nu_{u}$ whenever
$u\sim v$, we also have
$$
\sum_{(u,v)\in E(\gamma)}\nu_{uv}^{-1} \asymp \sum_{u\in
V(\gamma)}\nu_u
$$
Combining these inequalities together, we complete the proof. \qed\\

 Since
$\rho_x(r)$ is increasing in $r$, Lemma \ref{c:0822} immediately
implies Corollary \ref{c:0953w} as follows. Recall that
$B_{\rho}(x,r)$ is a $\rho-$ball. One can compare it with an
$L^1-$ball.
\begin{cor}\label{c:0953w} For  any $x\in \bZ^d$ and $r>0$,
$$
B(x, \rho^{-1}_x(c_1r))\subset B_{\rho}(x, r) \subset B(x,
\rho^{-1}_x(c_2r)).
$$
\end{cor}
Recall that $V_{\rho}(x,r)$ is the volume of $B_{\rho}(x,r)$. Set
$V(x,r)=\nu(B(x,r))$, similarly.

\begin{lemma}\label{l:0900} Let  $x\in \ZZ^d$ and $r>0$.
\begin{align*}
(1)&~V(x,r)\asymp r^d(|x|\vee r)^\al~~\hbox{if~}r\ge 1.\\
(2)&~V_{\rho}(x,r)\asymp V(x, \rho_x^{-1}( r))\asymp\begin{cases}\nu_x&~\hbox{if~~}r<\nu_x;\\
r^d|x|^{-(d-1)\al}
&~\hbox{if~~}\nu_x\le r\le |x|^{1+\al};\\
r^{(d+\al)/(1+\al)}&~\hbox{if~~}r> |x|^{1+\al}.
\end{cases}
\end{align*}
\end{lemma}
\pf  (1)  Let $x_1$ be the first coordinate of $x$ and set
$$ \Lambda=\{s=(s_1,\cdots,s_d)\in B(x,r):s_1=x_{1}\}. $$
 Write $e_1=(1,0,0,\cdots,0)\in \ZZ^d$.
 By Lemma  \ref{l:vip0926}, for each $s\in \Lambda$ we have
\begin{align*}
\sum_{l=-r}^{r}\nu_{s+l e_1}\asymp r(|s-re_1|\vee|s+re_1|)^\al\asymp
r(|s|\vee r)^\al\asymp r(|x|\vee r)^\al.
\end{align*}
Hence,
\begin{align}\label{e:1101}
V(x,r)=&\sum_{s\in \Lambda}\sum_{l=-r}^{r}\nu_{s+l e_1}\asymp
|\Lambda|\cdot r (|x|\vee r)^{\al} \asymp r^d(|x|\vee r)^\al.
\end{align}
(2) 
 Using (\ref{e:1101}) and Corollary \ref{c:0953w}, we get the
desired result.
 \qed\\

\begin{lemma}\label{l:1331} Let $w\in \bZ^d$ and  $R\ge 1$.  Then for any $x\in B(w,R)$
and $r\in [1, R]$,
\begin{equation}\label{e:0740ss}
V(w,R)\le c_1\left(\frac{R}{r}\right)^{c_1}V(x,r).
\end{equation}
Especially, $V(w,R)\le c_1 R^{c_1}\nu_x$.
\end{lemma}
\pf  It follows directly from Lemma \ref{l:0900} (1).
 \qed\\

So, $\nu(B(w,R))$ satisfy the volume doubling property in any case.
However, $\mu(B(w,R))=\sum_{x\in B(w,R)}\mu_{x}$ do not satisfy the
volume doubling property since $\mu(\bZ^d)<\infty$ when $\al> d$.\\

In \cite{VB}  Vir$\acute{a}$g, extending the early result of
\cite{DS}, showed that  Poincar$\acute{e}$ inequalities hold in any
convex lattices. We shall apply their technique to our weighted
lattices.
\begin{lemma}\label{T:LineP} Let  $x\in \ZZ^d$,  $r>0$.
Then for any function $f$ on  $B(x,r)$,
\begin{equation}\label{e:lineP}
 \min_a\sum_{u\in B(x,r)} (f(u)-a)^2 \nu_u\le  c_1 [\rho_x(r)]^2  \sum_{u,v\in  B(x,r)} (f(u)-f(v))^2
 \mu_{uv}.
\end{equation}
\end{lemma}
\pf   If $r\in(0,1)$ then $B(x,r)=\{x\}$ and  (\ref{e:lineP}) holds
since both side of the inequality are  zero. So, we may assume that
$r\ge 1$ in the following.

 By \cite[Proposition 2]{VB},   for each
$u,v\in \bZ^d$ we can choose a path $\gamma_{uv}$ such that, (1)
$\gamma_{uv}$ is an $L^1-$geodesic path from $u$ to $v$; (2) each
site in $\gamma_{uv}$ has $L^\infty-$distance less than 1 from the
Euclidean line $\overline{uv}.$ 
 For  $u,y\in \bZ^d$, write
$$ \Lambda_{uy}= \{s+z: s\in \gamma_{y,2y-u},
|z|\le 4, z\in \ZZ^d\}.
$$
By the construction, we have   $$ 1_{\{y\in \gamma_{uv}\}}\le
1_{\{v\in \Lambda_{uy}\}}+1_{\{u\in
\Lambda_{vy}\}}~{\rm~for~all~}u,v,y.
$$
By Lemma \ref{c:0822},
\begin{align*}
\sum_{v\in \Lambda_{uy} }  \nu_v\le& \sum_{s\in
\gamma_{y,2y-u}}\sum_{z\in \ZZ^d, |z|\le 4}\nu_{s+z}\le c \sum_{s\in
\gamma_{y,2y-u}}\nu_s\le c'\rho_y(|y-u|).
\end{align*}
So,  if  $u,y\in B(x,r)$,  we can use (\ref{e:1650}) and get
$$
\sum_{v\in \Lambda_{u,y} }  \nu_v\le c\rho_y(2r)\le c'\rho_x(r).
$$
By Lemma \ref{c:0822}, if $u,v\in B(x,r)$ then
\begin{equation}\label{e:1701ss}
\sum_{(y,z)\in E(\gamma_{uv})} \mu_{yz}^{-1}\le c\rho_u(|u-v|)\le
c\rho_x(r).
\end{equation}
Therefore, writing $B=B(x,r)$,
\begin{align*}
\sum_{u\in B} (f(u)-\overline{f})^2 \nu_u\le &
\frac{1}{\nu(B)}\sum_{u,v\in B} (f(u)-f(v))^2 \nu_u\nu_v
=\frac{1}{\nu(B)} \sum_{u,v\in B} \left(\sum_{(y,z)\in
E(\gamma_{uv})}
(f(y)-f(z))\right)^2 \nu_u\nu_v\\
\le& \frac{1}{\nu(B)} \sum_{u,v\in
B } \left(\sum_{(y,z)\in E(\gamma_{uv})} (f(y)-f(z))^2\mu_{yz}\right) \left(\sum_{(y,z)\in E(\gamma_{uv})} \mu_{yz}^{-1}\right)  \nu_u\nu_v \\
\le& \frac{c\rho_x(r)}{\nu(B)} \sum_{u,v\in
B }~\sum_{(y,z)\in E(\gamma_{uv})} (f(y)-f(z))^2\mu_{yz}  \nu_u\nu_v \\
\le&\frac{c\rho_x(r)}{\nu(B)} \sum_{y,z\in B} (f(y)-f(z))^2\mu_{yz}
\sum_{u,v\in B}1_{\{y\in \gamma_{u,v}\}}\nu_u\nu_v \\
 \le& \frac{c\rho_x(r)}{\nu(B)} \sum_{y,z\in B}
(f(y)-f(z))^2\mu_{yz} \left(\sum_{u\in B}\nu_u\sum_{v\in
\Lambda_{u,y} }  \nu_v+ \sum_{v\in B}\nu_v\sum_{u\in \Lambda_{v,y} }
\nu_u
\right) \\
\le& c'[\rho_x(r)]^2 \sum_{y,z\in B} (f(y)-f(z))^2\mu_{yz},
\end{align*}
 where  the second inequality is by the Cauchy-Schwarz inequality.  \qed

\begin{lemma} \label{L:DFest} Let $w\in \bZ^d$, $R\ge 1$ and $r\in (0, \rho_w(R)]$. Let
  $g:B(w,R)\to \bR^+$
with $\sum_{x\in B(w,R)} g(x) \nu_x \le 1$. Then \be
  \sum_{x,y\in B(w,R) }(g(x)-g(y))^2 \mu_{xy}
  \ge c_1r^{-2}  \left( \sum_{x\in B(w,R) }g(x)^2
  \nu_x-\frac{c_2}{V(w,R)}\left(\frac{\rho_w(R)}{r}  \right)^{c_3} \right).
\ee
\end{lemma}

\proof Let $\widehat{r}=\min\{\rho_x^{-1}(r): x\in B(w,R)\}\wedge
R$. Since $r\le \rho_w(R)$ and $\rho_x(R)\asymp \rho_w(R)$  for each
$x\in B(w,R)$, we have
\begin{equation}\label{e:0959ss}\rho_x(\widehat{r})\le c_1 r.\end{equation}
Note that for any $x\in B(w,R)$,
$$
\frac{\rho_x^{-1}(r)}{R}=\frac{\rho_x^{-1}(r)}{\rho_x^{-1}(\rho_x(R))}\ge
c_1\left(\frac{r}{\rho_x(R)}\right)^{c_1}\ge
c_2\left(\frac{r}{\rho_w(R)}\right)^{c_2}.
$$
So,   $\tfrac{\widehat{r}}{R}\ge
c_2\left(\frac{r}{\rho_w(R)}\right)^{c_2}$.  Using Lemma
\ref{l:1331}, we then have
\begin{equation}\label{e:0958s}
\frac{V(w,R)}{V(x,\widehat{r})}\le  c\left(\frac{R}{\widehat{r}}
\right)^c\le c_3\left(\frac{\rho_w(R)}{r}\right)^{c_3}.
\end{equation}

 Choose
$B_i=B(x_i,r_i), i=1,\cdots, N$  such that $B(w,R)=\cup_{i=1}^N
B(x_i, r_i) $ and
 $\widehat{r}\le r_i\le 2
\widehat{r}$ for each $i$, and
\begin{equation}\label{e:1001ss}
|\{i: x\in B(x_i,r_i)\}|\le c_4~~{\rm~for~all~}x\in B(w,R).
\end{equation}
Use Lemmas  \ref{T:LineP},
\begin{align*}
 \sum_{x,y\in\sB }(g(x)-g(y))^2 \mu_{xy}
  \ge& c_4^{-1}\sum_{i=1}^N  \sum_{x,y\in B_i }(g(x)-g(y))^2
  \mu_{xy}\\
  \ge& c_4^{-1}\sum_{i=1}^N  [\rho_{x_i}(\widehat{r})]^{-2}\sum_{x\in B_i }(g(x)-\overline{g}_i)^2
  \nu_x\\
   \ge& c_4^{-1} \sum_{i=1}^N (c_1r)^{-2} \left(\sum_{x\in B_i
   }g(x)^2\nu_x-\frac{(\sum_{x\in B_i} g(x)\nu_x)^2}{V(x_i,\widehat{r})}\right)\\
    \ge & (c_4c_1^2)^{-1}r^{-2}  \left( \sum_{x\in \sB }g(x)^2
  \nu_x-\frac{c_3}{V(w,R) }\left( \frac{\rho_w(R)}{r} \right)^{c_3}\sum_{i=1}^N \left(\sum_{x\in B_i}
  g(x)\nu_x\right)^2\right),
\end{align*}
where  $\sB=B(w,R)$ and  $\overline{g}_{i}$ is  the mean of $g$ on
$B_i$. Using (\ref{e:1001ss}), we get
$$
\sum_{i=1}^N \sum_{x\in B_i} g(x)\nu_x \le c\sum_{x\in
\sB}g(x)\nu_x\le c.
$$
Combining these inequalities with  $\sum_i a_i^2\le (\sum_i a_i)^2$
for all $a_i\ge 0$, we complete the proof.
 \qed\\

\begin{remark}\label{r:0817w} One cannot expect to improve
Lemma \ref{L:DFest} to the whole space such    as
 \begin{equation}\label{e:1255ss}
  \sum_{x,y\in \bZ^d }(g(x)-g(y))^2 \mu_{xy}
  \ge c_1r^{-2}  \left( \sum_{x\in \bZ^d }g(x)^2
  \nu_x-\frac{c_2}{V(w,R)}\left(\frac{\rho_w(R)}{r}  \right)^{c_3} \right)
\end{equation}
for all $r\in (0, \rho_w(R)]$,
 and  $g:\bZ^d\to \bR^+$
with $\sum_{x\in \bZ^d} g(x) \nu_x \le 1$.\\

 To see this, we fix
$\al\in (-1,0)$ and $d\ge 2$. On the one hand, choose $ R\ge 1$ and
$w\in \bZ^d$ with $|w|=R^{-\al^{-1}}$. Then $\rho_w(R)=1$, and hence
one can take $r=1$ further. Such,
\begin{equation}\label{e:1315ss}
V(w,R)\asymp R^{d-1}\rho_w(R)=R^{d-1}\rightarrow \infty.
\end{equation}
On the other hand, let $s\ge 1$,  and take
$$g(x)=A(s-|x|)1_{B(0,s)}(x),~~ x\in \bZ^d,$$ where $A$ is the
 constant which such that $\sum_x g(x)\nu_x=1$. Then
$$
 \sum_{x,y\in \bZ^d }(g(x)-g(y))^2 \mu_{xy}\le A^{2} \sum_{x,y\in B(0,s) }
 \mu_{xy}\le cA^{2}s^{d-\al},
$$
and
$$
\sum_{x\in \bZ^d }g(x)^2\nu_x \ge \frac{A^{2}s^2}{4}  \sum_{x\in
B(0,s/2)}\nu_x\ge cA^2s^{d+2+\al}.
$$
So, as $s$ goes to infinity,
\begin{equation}\label{e:1316ss}\sum_{x,y\in \bZ^d }(g(x)-g(y))^2
\mu_{xy}\ll \sum_{x\in \bZ^d }g(x)^2\nu_x.\end{equation} By
(\ref{e:1316ss}) and (\ref{e:1315ss}), the inequality
(\ref{e:1255ss}) fails.
\end{remark}

\section{On-diagonal upper bound estimates}
Fix $w\in \ZZ^d$, $R\ge 1$ and $T=\rho_w(R)^2$.    In this section,
our aim is to give an upper bound of  $p_T(w,w)$.  As Lemma
\ref{L:DFest} and Remark \ref{r:0817w} say,  we   have a good ball
$B(w,R)$ only. So, we turn to  the random walk $X$ with reflection
at $\partial_{i} B(w,R)$. By the approach of Barlow and Chen
\cite{BC}, we  obtain upper bounds  of the  heat kernel  of the
reflection process, and then bring these bounds back to the original
 process.\\

 Write $\sB=B(w,R)$ for short.  Let $Y$ be the continuous time random walk on $\sB$
with generator
$$
\mathscr{L}_{\sB} f(x)=\frac{1}{\nu_x}\sum_{y\in
\sB}(f(y)-f(x))\mu_{xy}.
$$
For $x\in  \ZZ^d$ and $r>0$, set
\begin{equation}\label{e:1525ss}
\tau_{x,r}=\inf\{t\ge 0: X_t\not\in B(x,r)\}.
\end{equation}
If $Y$ and $X$ start at the same vertex in $B(w,R-1)$, then  we can
couple $Y$ and $X$ on the same probability space such that
 \be
\label{e:XeqY10} Y_s=X_s~~~~ \hbox{ for ~ } 0\le s \le \tau_{w,R-1}.
\ee
 We use $\P_x$ for both $X$ and $Y$. Denote the heat kernel of $Y$
 by
$$
q_t(x,y)=\frac{\P_x(Y_t=y)}{\nu_y}.
$$

\begin{prop}\label{p:1612} For $u\in \sB$ and  $t\in (0,T]$,
\begin{align}
q_t(u,u)\le \label{e:2136}
\frac{c_1}{V(w,R)}\left(\frac{T}{t}\right)^{c_2}.
\end{align}
Especially, $q_T(w,w)\le \frac{c_1}{V(w,R)}.$
\end{prop}
\pf  Given Lemma \ref{L:DFest}, the proof is similar to
\cite[Proposition 3.1]{MB1}  and \cite[Proposition 3.2]{BC}
, so we
omit it. \qed\\

\begin{lemma}\label{L:1957}  Let
$x_1,x_2\in \sB$  with $|x_1-x_2|\ge \tfrac1{16} R$. If $t\le c_1T$
and $R\ge c_2$, then
\begin{equation}\label{e:1523w}
q_t(x_1,x_2)\le \tfrac1{4 V(w,R)}.
\end{equation}
\end{lemma}

\pf Write $\eta=\max_{x\in\sB}\nu_x$. By (\ref{e:0935}) and
(\ref{e:1650}), we have
\begin{align}\label{e:2045}
\frac{T^{1/2}}{\eta}=\frac{\rho_w(R)}{\max_{x\in
\sB}\rho_x(1)}=\inf_{x\in\sB}\left\{\frac{\rho_w(R)}{\rho_x(
R)}\cdot \frac{\rho_x(R)}{\rho_x(1)}\right\} \ge c_1 R^{c_1}.
\end{align}
Set $c_2=2^{|\al|+2}d$. Let $\widetilde{\nu}_x=\eta^{-1}\nu_x,
\widetilde{\mu}_{xy}=\eta \mu_{xy}$ and
$\widetilde{\rho}(x,y)=c_2^{-1}\eta^{-1}\rho(x,y)$ for $x,y\in \sB$.
Then
\begin{align}\label{e:0906w}\begin{cases}
\frac{1}{\widetilde{\nu}_x}\sum_{y\in\sB}\widetilde{\rho}(x,y)^2\widetilde{\mu}_{xy}\le
1;\\
~\\
 \widetilde{\rho}(x,y)\le 1 {\rm~~~whenever~}x\sim y.\end{cases}
\end{align}
Hence $\widetilde{\rho}(\cdot,\cdot)$ is an adapted metric, which
was introduced by Davies \cite{Da} and \cite{Da2}. Let
$Z_s=Y_{\eta^2s}$, for $s\ge 0$. Then $Z$ has the generator
$$
\widetilde{\mathscr{L}}_\sB
f(x)=\frac{1}{\widetilde{\nu}_x}\sum_{y\in \sB}
(f(y)-f(x))\widetilde{\mu}_{xy}.
$$
We state that  there exists constant $c,c'>0$ such that if $s\le
c\eta^{-2}T$ and $R\ge c'$ then
\begin{equation}\label{e:1624ss}
\P_{x_1}(Z_s=x_2)\le \frac{\nu_{x_2}}{4V(w,R)}.
\end{equation}
If this is true, then we have (\ref{e:1523w}) and prove the lemma.

We now prove  (\ref{e:1624ss}).
 Set
$c_3=c_{\ref{p:1612}.2}+c_1^{-1}c_{\ref{l:1331}.1}$. For each $i \in
\{1,2\}$, define
$$f_{x_i}(s)=\frac{V(w,R)}{c_{\ref{p:1612}.1}
\nu_{x_i}}\left(\frac{\eta^2s}{T}\right)^{c_3},~~~s\ge 0.$$ Then by
Proposition \ref{p:1612}, for  $s\le \eta^{-2}T$,
\begin{equation}\label{e:2003} \P_{x_i}(Z_s=x_i)= \P_{x_i}(Y_{\eta^2 s}=x_i)=q_{\eta^2s}(x_i,x_i)\nu_{x_i}\le\frac{1}{f_{x_i}(s)}.\end{equation}
Next we shall estimate the off-diagonal transition probability
$\P_{x_1}(Z_s=x_2)$ by using the 'two-point' method of
Grigor'yan-see
 \cite{Gr1,CGZ,FL,XC}.
The metric $d_\nu(x,y)$ in \cite{XC} is just $\widetilde{\rho}(x,y)$
and one can easily check that  $f_{x_i}(s)$ is $(1,2)$--regular on
$(0,T]$: see \cite{Gr1, XC} for the definition. By  (\ref{e:2045})
and  Lemma \ref{l:1331}, for $s\le \eta^{-2}T$,
\begin{align*} \frac{f_{x_i}(s)}{s^{c_3}}= &
\frac{V(w,R)}{c_{\ref{p:1612}.1}\nu_{x_i}}\cdot
\left(\frac{\eta^2}{T}\right)^{c_3}\le
\frac{c_{\ref{l:1331}.1}R^{c_{\ref{l:1331}.1}}}{c_{\ref{p:1612}.1}}\cdot (c_1R^{c_1})^{-2c_3}\\
=& c' R^{c_{\ref{l:1331}.1}-2c_1c_3}\le c'R^{-c_{\ref{l:1331}.1}}\le
c'.
\end{align*}
 Therefore, by \cite[Theorem 1.1]{XC}  for  $s\in \big(~
\widetilde{\rho}(x_1,x_2),~\eta^{-2} T\big]$,
\begin{align}
\label{e:0947ss} \P_{x_1}(Z_s=x_2) \le& \frac{c_4
(\widetilde{\nu}_{x_2}/\widetilde{\nu}_{x_1})^{1/2}} {\sqrt{f_{x_1}(
c_5 s )f_{x_2}( c_5
s)}}\exp\left(-c_6\frac{{\widetilde{\rho}}(x_1,x_2)^2}{s}\right)\\
\label{e:0949} =&\frac{c_7
\nu_{x_2}}{V(w,R)}\left(\frac{T}{\eta^2s}\right)^{c_3} \exp\left(
-c_6c_2^{-2}\frac{{\rho}(x_1,x_2)^2}{\eta^2s}\right).
\end{align}
By Lemma \ref{c:0822} and the condition  $|x_1-x_2|\ge
\tfrac1{16}R$, we have
\begin{equation}\label{e:14481w}\rho(x_1,x_2)\ge
c_{\ref{c:0822}.2}\rho_{x_1}(\tfrac1{16}R)\ge
c_8\rho_w(R)=c_8T^{1/2}.\end{equation} Substituting (\ref{e:14481w})
into (\ref{e:0949}) gives
$$
\P_{x_1}(Z_s=x_2) \le \frac{c_7
\nu_{x_2}}{V(w,R)}\left(\frac{T}{\eta^2s}\right)^{c_3} \exp\left(
-c_6c_2^{-2}c_8^2\frac{T}{\eta^2s}\right),
$$
which implies (\ref{e:1624ss}) holds for each $s\in \big(~
\widetilde{\rho}(x_1,x_2),~c_{9}\eta^{-2} T\big]$, provided
$c_{9}>0$ is small enough.

On the other hand, by \cite[Corollary 2.8]{XC} we have the `long
range' bounds, that is,  if $s\le \widetilde{\rho}(x_1,x_2)$ then
\begin{equation}\label{e:1446w}
\P_{x_1}(Z_{s}=x_2) \le c' (\widetilde{\nu}_{x_2}/
\widetilde{\nu}_{x_1})^{1/2} e^{-c\widetilde{\rho}(x_1,x_2)}.
\end{equation}
Using  (\ref{e:2045}) and (\ref{e:14481w}), we have
\begin{equation}\label{e:1525w}
\widetilde{\rho}(x_1,x_2)=c_2^{-1}\eta^{-1}\rho(x_1,x_2)\ge
c\eta^{-1} T^{1/2}\ge c'R^{c''}.\end{equation} Combining these
inequalities with Lemma \ref{l:1331},
\begin{align*}
\P_{x_1}(Z_{s}=x_2) \le& c ({\nu}_{x_2}/ {\nu}_{x_1})^{1/2}
e^{-c'R^{c''}}= \frac{c\nu_{x_2}}{V(w,R)}\cdot
\frac{V(w,R)}{(\nu_{x_2}\nu_{x_1})^{1/2}}e^{-c'R^{c''}}\\
\le& \frac{c\nu_{x_2}}{V(w,R)}\cdot
c_{\ref{l:1331}.1}R^{c_{\ref{l:1331}.1}}e^{-c'R^{c''}}.
\end{align*}
So,  (\ref{e:1624ss}) holds  again if $s\le
\widetilde{\rho}(x_1,x_2)$ and $R\ge c$.
 \qed\\

\begin{lemma}\label{L:heatZ10}
Let $  t\le  c_1T$ and $x\in B(w,\tfrac{7}{8}R)$. If $R\ge c_2$ then
$$ \P_x \big(Y_t\not\in B(x,\tfrac1{16}R)\big) \le \tfrac14. $$
\end{lemma}

\pf By Lemma \ref{L:1957}, we get
$$
 \bP_x\big( Y_t \not\in  B(x,\tfrac1{16}R) \big)=\sum_{y \in \sB-
 B(x,\tfrac1{16}R)}q_t(x,y)\nu_y
\le \sum_{ y \in \sB-   B(x,\tfrac1{16}R)}
 \frac{\nu_y}{4V(w,R)}
\le \frac{1}{4}.
$$
 \qed\\

 Now we bring
these bounds of the reflection process back to the original process.
Note that $X$ and $Y$ agree until time $\tau_{w,R-1}$.

\begin{lemma}\label{l:1020}
 If $R\ge c_1$ then for  $x\in B(w, \tfrac58R)$,
$$ \P_x(\tau_{x, R/8}<    c_2T)\le  \tfrac{1}{2}. $$
\end{lemma}

\pf  Given Lemma \ref{L:heatZ10}, the proof is similar to
\cite[Lemma 4.1]{BC}, so we
omit it. \qed\\

\begin{prop}\label{p:0812} Let $w\in \ZZ^d$, $R>0$ and
$T=\rho_w(R)^2$. Then $$\P_w(X_{T}=w)\le \frac{c_1\nu_w}{V(w,R)}.$$
\end{prop}
\pf If $R<(c_{\ref{l:1020}.1}\vee c_{\ref{L:1957}.2})$ then by Lemma
\ref{l:1331},
$$
\frac{\nu_w}{V(w,R)}\ge c'R^{-c}\ge c'(c_{\ref{l:1020}.1}\vee
c_{\ref{L:1957}.2})^{-c} \ge c_1^{-1}\P_w(X_{T}=w).
$$
So, let $R\ge (c_{\ref{l:1020}.1}\vee c_{\ref{L:1957}.2})$. Given
Lemma \ref{l:1020}, similar to the inequality  (4.6) of Barlow and
Chen \cite{BC} we obtain
$$
p_{ c_2 T}(w,w)\le q_{c_2 T}(w,w) + \sup_{0< s \le c_2 T}\max_{y\in
A}q_s(y,w),
$$
where $c_2= c_{\ref{l:1020}.1}\wedge c_{\ref{L:1957}.1}\wedge1$ and
$ A= B(w, 5R/8)-B(w,5R/8-1).$ By Proposition \ref{p:1612} and Lemma
\ref{L:1957},
$$
p_T(w,w)\le p_{ c_2 T}(w,w)\le \frac{c_3}{V(w,R)}.
$$
 \qed
\section{Near diagonal lower bound estimates}
 In this section, we shall prove the
following lower bounds for the near diagonal transition
probabilities. Recall $\tau_{x,r}$ from section 3. Fix $\delta\in
(0,1/2)$. We will use the notation $K_i$ to denote constants which
depend only $\delta,\al$ and $d$, while
 $c_i=c_i(\al,d)$ as before.

\begin{thm}\label{t:09051}
   Let $w\in \bR^d$ and   $R\ge 1$. For
$x_1, x_2\in B(w,R)$ and $t\in [\delta \rho_w(R)^2, ~2\rho_w(R)^2]$,
\begin{equation}\label{e:2036q}
\P_{x_1}(X_t=x_2, ~\tau_{w, c_1R}>t)\ge K_2\frac{\nu_{x_2}}{V(w,R)}.
\end{equation}
\end{thm}

Since $\mu(B(w,R))$ do not satisfy the volume doubling property, we
cannot obtain the lower bound by a general approach.  Let us begin
with  a ball   far from the origin.

\begin{lemma}\label{e:21102}   Let  $w\in\bZ^d$ and $R\ge 1$ with $|w|\ge 32R$. Then  for any $x_1, x_2\in B(w,R)$ and
$t\in [\delta  \rho_w(R)^2, 2\rho_w(R)^2]$,
\begin{equation}\label{e:2036}
\P_{x_1}(X_t=x_2, ~\tau_{w,8R}>t)\ge K_1 R^{-d}.
\end{equation}
\end{lemma}
\pf  Since $|w|\ge 32 R$, $\rho_w(R)=R|w|^\al$, moreover, for any
$x,y\in B(w,16R)$ with $x\sim y$,
$$
\nu_x\in [4d c_1^{-1} |w|^\al, c_1
|w|^{\al}]~~~{\rm~and~~~}\mu_{xy}\in [c_1^{-1}|w|^{-\al}, c_1
|w|^{-\al}].
$$
By the application of  Lemma \ref{l:1020} on $B(w, 8R)$,  there
exists  $c_2\in (0,1/2)$ such that
\begin{equation}\label{e:1010ss}
\P_{x}(\tau_{x, R}>  c_2\rho_w(R)^2)\ge
\tfrac{1}{2},~~{\rm~for~all~}x\in B(w,R).
\end{equation}
For each $x,y\in B(w, 16R)$, we set
\begin{align*}
\widetilde{\nu}_x=c_1|w|^{-\al}\nu_x~~~~{\rm and~~~~}\widetilde{\mu}_{xy}=\begin{cases}c_1^{-1} |w|^\al\mu_{xy},& {\rm~if~}x\not= y;\\
\widetilde{\nu}_x-c_1^{-1} |w|^\al\sum_{z\in  B(w,16R)\setminus\{
x\}}~\mu_{xz},& {\rm~if~}x=y.
\end{cases}
\end{align*}
So,  $ \widetilde{\nu}_x,\widetilde{\mu}_{xy} \in [c_3, c_3^{-1}] $
for all $x\in B(w,16R)$ and $y\in  B(w,16R)\cap B(x,1)$.
  Let
$Z$ be the continuous time (constant speed) random walk on $
B(w,16R)$ with generator
$$
\widetilde{\mathscr{L}} f(u)=\frac{1}{\widetilde{\nu}_u}\sum_{v\in
 B(w,16R)} (f(v)-f(u))\widetilde{\mu}_{uv}.
$$
Then  $Z$ and $X$ can be coupled in the same probability such that
$$Z_{s}=X_{c_1^{-2}|w|^{2\al} s},~~~{\rm~for~all~}s<\sigma=c_1^{2}|w|^{-2\al} \tau_{w,8R},$$
where $\sigma:=\inf\{s\ge 0: Z_s\not\in B(w,8R)\}.$ Fix $x_1,x_2\in
 B(w,R)$, and let $ u(s,y)=\P_{x_1}(Z_s=y,\sigma>s)/
\widetilde{\nu}_y$ for each $y\in  B(w,16R)$ and $s\ge 0$. Then $u$
is a positive solution of the heat equation $\frac{\partial
u}{\partial s}=\widetilde{\mathscr{L}} u$ on $(0,\infty)\times
 B(w,4R)$. One can easily check that $DV(C_1), P(C_2)$ and
$\Delta(\al)$ hold for the weighted graph with vertex set $B(w,16R)$
and edge weight $\widetilde{\mu}_{xy}$, and so $u(s,y)$ satisfies
the Harnack inequality, see \cite[Theorem 1.7]{DT}. Therefore,
$$
 \max_{[ \tfrac{1}{2}s_0,  s_0]\times B(w,2R)} u \le K_1^{-1}~ \min_{[\delta  c_1^2R^2, 2c_1^2R^2]\times
B(w,2R)} u,
$$
where $s_0=\delta c_2 c_1^2 R^2$.
 Furthermore, for any $s\in [\delta   c_1^2 R^2,
2c_1^2R^2]$,
\begin{align}
\nonumber\P_{x_1}(Z_s=x_2,\sigma>s)\ge&K_1\left(\sum_{z\in
B(w,2R)}\widetilde{\nu}_z \right)^{-1}\sum_{z\in
B(w,2R)}\P_{x_1}(Z_{s_0}=x_2,\sigma>s_0)\\
\nonumber\ge&K_1 c_3 |B(w,2R)|^{-1}\P_{x_1}(Z_{s_0}\in B(w,2R),\sigma>s_0)\\
\label{e:1434ss}\ge& K_1 c_3(5R)^{-d}\P_{x_1}(\inf\{h:Z_h\not\in
B(x_1, R)\}> s_0).
\end{align}
Since  $X_{ t}=Z_{c_1^2|w|^{-2\al}t}$ for all $t<\tau$, inequality
(\ref{e:1434ss}) can be rewrote as
$$
\P_{x_1}(X_t=x_2,\tau_{w,8R}>t)\ge K_2R^{-d}\P_{x_1}(\tau_{x_1,R}>
\delta c_2  R^2|w|^{2\al}),~~t\in [ \delta  R^{2}|w|^{2\al},
2R^{2}|w|^{2\al}].
$$
Using (\ref{e:1010ss}), we  finish the proof.\qed\\

\begin{lemma}\label{l:1503ss} For any $t\in [\delta
R^{2+2\al}, R^{2+2\al}]$ and $x\in B(0,R)$,
$$
\P_x(|X_t|>K_1R, ~\tau_{x, c_2R}>t)\ge \tfrac14.
$$
\end{lemma}
\pf   By    Proposition \ref{p:0812}, for any  $x,y\in \bZ$ and
$t>0$,
\begin{align}
\nonumber p_t(x,y)\le (p_t(x,x)p_t(y,y))^{1/2}\le&
c_1(V(x,\rho_x^{-1}(t^{1/2}))V(y,\rho^{-1}_y(t^{1/2})))^{-1/2}.
\end{align}
So, from  Lemma \ref{l:0900} we can get, if $x,y\in B(0,  R)$ and
$t\in [\delta  R^{2+2\al}, R^{2+2\al}]$ then
$$
p_t(x,y)\le K_1t^{-(d+\al)/(2+2\al)}\le K_2R^{-d-\al}.
$$
Fix  $x\in B(0,  R)$ and $t\in [\delta  R^{2+2\al}, R^{2+2\al}]$. By
Lemma \ref{l:0900} again, for each $\epsilon\in(0,1)$,
\begin{align*}
\P_x(|X_t|\le \epsilon R)=&\sum_{y\in B(0,\epsilon R)}
p_t(x,y)\nu_y\le  V(0,\epsilon R)\cdot K_2R^{-(d+\al)}\\
\le& c_2(\epsilon R)^{d+\al}\cdot
K_2R^{-(d+\al)}=K_2c_2\epsilon^{d+\al}.
\end{align*}
Hence there exists $\epsilon_0=\epsilon_0(\delta,\al,d)>0$ such that
\begin{equation}\label{e:1421ss}
\P_x(|X_t|\le \epsilon_0 R)\le \tfrac{1}{4}.
\end{equation}

On the other hand, applying  Lemma \ref{l:1020} gives
\begin{equation}\label{e:1904w}
\P_x(\tau_{x, c R}<  t)\le \P_x(\tau_{x, c R}< R^{2+2\al})\le
\tfrac12.
\end{equation}
Combing  (\ref{e:1904w}) with (\ref{e:1421ss}), we finish the proof.
\qed

\begin{lemma}\label{e:21102ss}   Let  $R\ge 1$. Let $x_1, x_2\in B(0,R)\setminus B(0,\delta R)$ and $t\in [\delta  R^{2+2\al},
R^{2+2\al}]$. Then
\begin{equation}\label{e:2036}
\P_{x_1}(X_t=x_2, ~\tau_{0,10R}>t)\ge K_1 R^{-d}.
\end{equation}
\end{lemma}
\pf  Write $\mathbb{T}=B(0,R)\setminus B(0,\delta R)$ for short.  If
$d\ge 2$, then $\mathbb{T}$ is connected. Note that
$\rho_w(R)\subset [K_1^{-1}R^{1+\al}, K_1R^{1+\al}]$ and
$B(w,8R)\subset B(0,10R)$ for all $w\in \mathbb{T}$, and there exist
vertices $ w_i\in \mathbb{T},~i\le K_2$ such that
$\mathbb{T}=\cup_{i=1}^{K_2} B(w_i, \tfrac{\delta}{64}R)$.  A
standard chaining argument using Lemma \ref{e:21102} on
$B(w_i,\tfrac{\delta}{32}R)$, proves (\ref{e:2036}) for $d\ge 2$.
Next, we consider  $d=1$. Since  $\mathbb{T}=([-R,-\delta
R]\cup[\delta R, R])\cap \bZ$ is not connected, we  have to discuss
the problem on several cases.

 Case I: $x_1,x_2> 0$.  Then $x_1$ and $x_2$ can be joint with a
sequence of balls $B(w_i,\tfrac{\delta}{32}R)$ within $[\delta R,
R]\cap \bZ$ as before. Hence  (\ref{e:2036}) holds for this case,
too.

  Case II: $x_1>0>x_2$. For conciseness, we write
$\widehat{\P}$ for the measure of the process $X$ killed on exiting
$B(0, 10R)$. Let $\epsilon_0=\epsilon_0(\delta,\al,d)\in(0,1)$ be a
small constant, whose value will be taken later. Set $x_*=\lfloor
\epsilon_0 R\rfloor$.
 By the result of Case I, we
have
$$
\inf_{s\in [\delta'R^{2+2\al}, R^{2+2\al}]}~\inf_{x,y\in
B(0,R)\setminus B(0, \delta')}\widehat{\P}_{x}(X_s=y)\ge K_3 R^{-1},
$$
where $\delta'=\min\{\tfrac12\epsilon_0, \tfrac13\delta\}$. So,
$$
\widehat{\P}_{x_1}(X_{t/3}\in(\tfrac12\epsilon_0R,\epsilon_0 R))\ge
\tfrac14K_3\epsilon_0~~{\rm~and~~}\inf_{s\in [ t/3,
t]}\widehat{\P}_{-x_*}(X_{s}=x_2)\ge K_3 R^{-1}.
$$
For  $x\in \bZ$, we define $ \sigma_x=\inf\{t\ge 0: X_t=x\}$, the
 first time of  visiting vertex $x$. By the strong Markov property,
\begin{align}
\nonumber\P_{x_1}(X_t=x_2, ~\tau_{0,10R}>t)\ge&
\widehat{\P}_{x_1}(\sigma_{x_*}<\tfrac{t}{3},
\sigma_{-x_*}<\tfrac{2t}{3},
X_{t}=x_2)\\
\nonumber\ge&\widehat{\P}_{x_1}(\sigma_{x_*}<\tfrac{t}{3})\widehat{\P}_{x_*}(
\sigma_{-x_*}<\tfrac{t}{3})\inf_{s\in [ t/3,
t]}\widehat{\P}_{-x_*}(X_{s}=x_2)\\
\nonumber\ge&\widehat{\P}_{x_1}(X_{t/3}\in(\tfrac12\epsilon_0R,\epsilon_0
R))\widehat{\P}_{x_*}( \sigma_{-x_*}<\tfrac{t}{3})\inf_{s\in [ t/3,
t]}\widehat{\P}_{-x_*}(X_{s}=x_2)\\
\label{e:1849ss}\ge&\tfrac14K_3\epsilon_0 \cdot\widehat{\P}_{x_*}(
\sigma_{-x_*}<\tfrac{t}{3})\cdot K_3R^{-1}.
\end{align}
So, we need a lower bound of  $\widehat{\P}_{x_*}(
\sigma_{-x_*}<\tfrac{t}{3})$.  By  Lemma \ref{c:0822}, for any $x\in
\bN$ and $r,s\ge 2|x|$,
\begin{equation}\label{e:1104ss}\P_x(\sigma_{x-r}>\sigma_{x+s})=\frac{\sum_{i=x-r}^{x-1} \mu_{i,i+1}^{-1}}{\sum_{i=x-r}^{x+s-1} \mu_{i,i+1}^{-1}}\le
\frac{c_{\ref{c:0822}.1}^{-1}\rho_{x}(r)}{c_{\ref{c:0822}.1}\rho_{x}(r+s-1)}
 \le c_1\left(\frac{r}{r+s}\right)^\al.
\end{equation}
So, there exists $c_2\in \bN$ such that
\begin{equation}\label{e:1527ss}
\P_x(\sigma_{-x}>\sigma_{c_2x})\le \tfrac{1}{8},~~{\rm for~all~}x\in
\bN.
\end{equation}
By  Lemma \ref{l:1503ss}, there exist $c_3=c_3(\al,d)\in(0,1)$ and
$K_4=K_4(\delta,\al,d)\in (0,1)$ such that
$$
\P_x(|X_{t/3}|>K_4 R, ~~\tau_{x, R}>t/3)\ge \tfrac14,~~x\in B(0,
c_3R).
$$
  Now we choose $\epsilon_0=c_2^{-1}K_4c_3$.    Then $x_*=\lfloor c_2^{-1}K_4c_3R\rfloor
\in B(0, c_3R)$ and so,
\begin{equation}\label{e:1526ss}
{\P}_{x_*}(\sigma_{-x_*}\wedge\sigma_{c_2x_*}<\tfrac13t,~\tau_{0,10R}>\tfrac13t)\ge
\P_{x_*}(|X_{t/3}|> K_4R,~ \tau_{c_2x_*, R  }>t/3)\ge \tfrac14.
\end{equation}
Combining (\ref{e:1526ss}) with (\ref{e:1527ss}), we get
$$
\widehat{\P}_{x_*}(\sigma_{-x_*}<\tfrac t3)\ge \tfrac{1}{8}.
$$
Substituting the above inequality into (\ref{e:1849ss}),  we prove
(\ref{e:2036}) for the second case.

By symmetry, we  have (\ref{e:2036}) as $x_1<0$. Therefore,
(\ref{e:2036})  holds in any case.
 \qed\\

{\it Proof of Theorem \ref{t:09051}. } If  $|w|\ge 32R$, then one
can take $c_1=8$ in (\ref{e:2036q}) and the problem is reduced to
Lemma \ref{e:21102}. So, let $R>|w|/32$ in the following. Then
$$
\rho_w(R)\in [c_1 R^{1+\al},
c_2R^{1+\al}]~~{\rm~and~~~}B(w,R)\subset B(0, 40R).
$$
Fix $t\in [\delta\rho_w(R)^2,2\rho_w(R)^2]$. Then $t\in [c_1\delta
R^{2+2\al}, c_2 R^{2+2\al}]$. By Lemma \ref{l:1503ss},
 for any $x\in B(0,40R)$,
\begin{equation}\label{e:2218ss}
\P_x(|X_{t/3}|>K_1R,~ \tau_{0,c_3R}>t/3)\ge \tfrac14.
\end{equation}
 Write $\mathbb{T}=B(0, c_3R)\setminus B(0, K_1 R)$.
 By Lemma \ref{e:21102ss},  for all $x,y\in
\mathbb{T}$,
$$
\P_x(X_{t/3}=y,~\tau_{0,10c_3R}>t/3)\ge K_2 R^{-d}.
$$
 Therefore, for  any $x_1,x_2\in B(w,R)\subset B(0,
40R)$,
\begin{align}
\nonumber\P_{x_1}(X_t=x_2, \tau_{0, 10c_3R}>t) \ge& \sum_{x,y\in
\mathbb{T}} \widehat{\P}_{x_1}(
X_{t/3}=x,~ X_{2t/3}=y, ~X_t=x_2)\\
\nonumber=& \sum_{x,y\in \mathbb{T}}\widehat{ \P}_{x_1}(
X_{t/3}=x)\widehat{\P}_{x}(
X_{t/3}=y)\widehat{\P}_{y}( X_{t/3}=x_2)\\
\nonumber\ge& K_2R^{-d}\sum_{x,y\in \mathbb{T}} \widehat{\P}_{x_1}(
X_{t/3}=x)\widehat{\P}_{y}( X_{t/3}=x_2)\\
\nonumber=&K_2R^{-d}\sum_{x,y\in \mathbb{T}} \widehat{\P}_{x_1}(
X_{t/3}=x)\widehat{\P}_{x_2}( X_{t/3}=y)\frac{\nu_{x_2}}{\nu_y}\\
\nonumber\ge & K_2R^{-d}\frac{\nu_{x_2}}{\max_{y\in
\mathbb{T}}\nu_y}\sum_{x,y\in \mathbb{T}} \widehat{\P}_{x_1}(
X_{t/3}=x) \widehat{\P}_{x_2}( X_{t/3}=y)\\
\nonumber\ge & K_2R^{-d}\frac{\nu_{x_2}}{K_3R^{\al}}\widehat{ \P}_{x_1}(|X_{t/3}|\in \mathbb{T})\widehat{\P}_{x_2}(|X_{t/3}|\in \mathbb{T})\\
\label{e:1714q}\ge& \frac{K_2}{16K_3} \frac{\nu_{x_2}}{R^{d+\al}},
\end{align}
where we  use $\widehat{\P}$ to denote the measure of the process
$X$ killed on exiting $B(0, 10c_3R)$.
 Substituting $V(w,R)\le c R^{d+\al}$ and $\tau_{0,10c_3R}\le \tau_{w,
cR}$ into (\ref{e:1714q}), we  complete the proof.\qed\\

\section{Proof of Theorem \ref{T:lower1848}}

\begin{lemma}\label{l:1820w} There exists
constant $c_1>0$ such that  for any $x,y\in \bZ^d$,
$$
(\nu_x\vee \nu_y)|\log \nu_x-\log \nu_y|^3\le c_1 \rho(x,y).
$$
\end{lemma}
\pf Let $|x|> |y|\ge 1$. Directly calculate
\begin{align*}
\frac{(\nu_x\vee \nu_y)}{(|x|\vee |y|)^\al}\cdot\frac{|\log
\nu_x-\log \nu_y|^3}{ |x-y|}
=&\frac{|x|^\al\vee |y|^\al}{|x|^\al}\cdot\frac{|\log (|x|^\al)-\log (|y|^\al)|^3}{|x-y|}\\
=& (\tfrac{|y|}{|x|})^{\al\wedge0}\cdot |\al|^3 \frac{\log^3
(|x|/|y|)}{(|x|/|y|-1)|y|}\\
\le& |\al|^3 \sup_{t> 1}\left\{ t^{(-\al)\vee0}\cdot\frac{\log^3
t}{t-1}\right\}.
\end{align*}
Since $\al>-1$, the supremum  of the right side is finite and hence
if $|x|> |y|\ge 1$ then
$$
(\nu_x\vee \nu_y)|\log \nu_x-\log \nu_y|^3\le c(|x|\vee
|y|)^\al|x-y|\le c'\rho(x,y).
$$

 The proof of the rest case  is the same and so we
omit the
details. \qed\\
\\

{\it Proof of Theorem \ref{T:lower1848}.} We obtain the  Gaussian
upper bounds by the same way as  Lemma \ref{L:1957}. Write
$\eta=\nu_x\vee \nu_y$ for short. Set
$\widetilde{\nu}_u=\eta^{-1}\nu_v$ and $\widetilde{\mu}_{uv}=\eta
\mu_{uv}$   for each $u,v\in \ZZ^d$. Denote
$\widetilde{\rho}:\ZZ^d\times\ZZ^d\rightarrow \bR^+$ by
$$
\widetilde{\rho}(u,v)=\big((2^{|\al|+2}d)^{-1}\cdot\eta^{-1}\rho(u,v)\big)\wedge
|u-v|_1.
$$
Then $\widetilde{\rho}(\cdot,\cdot)$ is an adapted metric of
$\bZ^d$, that is, for all $u\in \bZ^d$,
\begin{align*}\begin{cases}
\frac{1}{\widetilde{\nu}_u}\sum_{v\in\bZ^d}\widetilde{\rho}(u,v)^2\widetilde{\mu}_{uv}\le
1;\\
~\\
 \widetilde{\rho}(u,v)\le 1 {\rm~~~whenever~}v\sim u.\end{cases}
\end{align*}
 By  Lemma \ref{c:0822}, $$ \eta^{-1}\rho(x,y)\le
(\nu_x\vee \nu_y)^{-1} \cdot c (|x|\vee |y|)^\al |x-y|\le
c|x-y|_1.$$ So,
\begin{equation}\label{e:94302}
c_1^{-1} \eta^{-1} \rho(x,y)\le \widetilde{\rho}(x,y)\le c_1
\eta^{-1} \rho(x,y).
\end{equation}
Set  $Z_s=X_{\eta^2s}$ for $s\ge 0$. Then $Z$ has the generator
$$
\widetilde{\mathscr{L}} f(u)=\frac{1}{\widetilde{\nu}_u}\sum_{v\in
\bZ^d} (f(v)-f(u))\widetilde{\mu}_{uv}.
$$
 By Proposition
\ref{p:0812} for each $z\in \{x,y\}$,
$$
\P_{z}(Z_s=z)=\P_{z}(X_{\eta^2s}=z)\le
\frac{c_2\nu_{z}}{V(z,\rho_{z}^{-1}(\eta
s^{1/2}))}:=\frac{1}{f_z(s)}.
$$
By  Lemma \ref{l:1331} and the inequality (\ref{e:0935}), for each
$s\ge (\log \nu_x-\log \nu_y)^2$ we have
\begin{align*}
f_z(s)\le& \frac{V(z,\rho_{z}^{-1}(\eta
e^{s^{1/2}}))}{c_2\nu_{z}}\le c(\rho_{z}^{-1}(\eta
e^{s^{1/2}}))^c\le c'\left(\frac{\eta e^{s^{1/2}}}{\rho_z(1)}
\right)^{c'}\le c'\left(\frac{\nu_x\vee \nu_y}{\nu_x\wedge \nu_y}
e^{s^{1/2}}\right)^{c'} \le c'e^{2c's^{1/2}}.
\end{align*}
Therefore, similar to (\ref{e:0947ss}) we can apply \cite[Theorem
5.1]{XC} and get
$$
\P_{x}(Z_{s}=y)\le  \frac{c_3
(\nu_{y}/\nu_x)^{1/2}}{\sqrt{f_z(s/c_3)f_y(s/c_3)}} \exp\left(
-\frac{\widetilde{\rho}(x,y)^2}{c_3s}\right)~{\rm~for~all~}s\ge
|c_3\log (\nu_x/\nu_y)|^3\vee \widetilde{\rho}(x,y).
$$
 By   the inequality (\ref{e:94302}) and Lemma \ref{l:1820w},
$$|c_3\log (\nu_x/\nu_y)|^3\le c_4\widetilde{\rho}(x,y).$$
So, for each $t\ge  c_4c_1\eta \rho(x,y)$, we have $\eta^{-2}t\ge
 |c_3\log (\nu_x/\nu_y)|^3\vee \widetilde{\rho}(x,y)$ and
\begin{align}
\nonumber\P_{x}(X_{t}=y)=&\P_{x}(Z_{\eta^{-2}t}=y)\le \frac{c_5
(\nu_{y}/\nu_{x})^{1/2}}{\sqrt{f_x(\eta^{-2} t/c_5)
f_y(\eta^{-2}t/c_5)}} \exp\left(
-\frac{\widetilde{\rho}(x,y)^2}{c_5\eta^{-2}t}\right)\\
\label{e:0931w}\le &\frac{c \nu_{y}}{\sqrt{V(x,\rho_{x}^{-1}(
t^{1/2}))V(y,\rho_{y}^{-1}( t^{1/2}))}} \exp\left(
-c'\frac{\rho(x,y)^2}{t}\right).
\end{align}
Further, by Lemma \ref{l:0900} we conclude  that
\begin{equation}\label{e:0931ss}
p_t(x,y)\le \frac{c}{\sqrt{V_{\rho}(x, t^{1/2})V_{\rho}(y, t^{1/2})}
}\exp\left(-c'\frac{\rho(x,y)^2}{t}\right), ~~t\ge c_4c_1\eta
\rho(x,y).
\end{equation}

On the other hand, by \cite[Corollary 2.8]{XC}, if $s\le
c_4c_1^2\widetilde{\rho}(x,y)$ then
\begin{align*}
\P_{x}(Z_s=y)\le&
c(\widetilde{\nu}_{y}/\widetilde{\nu}_x)^{1/2}\exp\Big(-c'\widetilde{\rho}(x,y)\big(1\vee
\log\left(\widetilde{\rho}(x,y)/s\right)\big)\Big).
\end{align*}
Hence, for each $t\le c_4c_1\eta \rho(x,y)$,
\begin{align}
 \nonumber\P_{x}(X_t=y)=&\P_x(Z_{\eta^{-2}t}=y)\le c(\widetilde{\nu}_{y}/\widetilde{\nu}_x)^{1/2}\exp\Big(-c'\widetilde{\rho}(x,y)\big(1\vee
\log\left(\eta^2\widetilde{\rho}(x,y)/t\right)\big)\Big)\\
\label{e:1111ss}\le&
c({\nu}_{y}/{\nu}_x)^{1/2}\exp\Big(-c''\eta^{-1} \rho(x,y)\big(1\vee
\log\left( \eta \rho(x,y)/ t\right)\big)\Big).
\end{align}
Combining (\ref{e:1111ss}) with (\ref{e:0931ss}), we conclude that
both (\ref{e:1843w}) and (\ref{e:1845w}) are true.\\

The Gaussian lower bound is proved by a standard chaining argument.
 If $t\ge \rho(x,y)^2$, then there exists $c_1>1$ such that
 $t\ge c_1^{-2}\rho_x(|x-y|)^2$.  Applying Theorem
\ref{t:09051} on $B(x, \rho_x^{-1}(c_1t^{1/2}))$, we get
\begin{equation}\label{e:1236wq}
p_t(x,y)\ge \frac{c}{V(x,\rho_x^{-1}(c_1t^{1/2}))}\ge
\frac{c'}{V_{\rho}(x,t^{1/2})}.
\end{equation}
So, let $(\nu_x\vee \nu_y)\rho(x,y)\le t\le \rho(x,y)^2$.  Fix an
$L_1-$geodesic path  $\gamma$ from $x$ to $y$.  By Lemma
\ref{c:0822}, there exists $c_2>1$ such that
$$
\nu(\gamma)\le c_2\rho(x,y).
$$
Set $r= {t}/{\rho(x,y)}$,
 then
$$
\rho(x,y)\ge r\ge  \nu_x\vee \nu_y=\max_{z\in \gamma}\nu_z.
$$
Hence  there exists a sequence of vertices $y=z_0,z_1\cdots,z_m=x$
on the path  $\gamma$, such that
$$
m\le 2c_2\rho(x,y)/r=2c_2\frac{\rho(x,y)^2}{t}~~~{\rm and~~~~}r\le
\rho(z_{i-1}, z_i)\le 2r~~{\rm~for~}i\le m.
$$
As a result,
$$
 |z_{i-1}-z_{i-2}|\le
c_3'\rho_{z_{i-1}}^{-1}(\rho(z_{i-1},z_{i-2}))\le
c_3'\rho_{z_{i-1}}^{-1}(2\rho(z_{i-1},z_{i}))\le
(c_3-1)|z_i-z_{i-1}|.
$$
Write $r_i=|z_i-z_{i-1}|$, $F_i=B(z_i,r_i)$ and $F_i^*=B (z_i,
c_3r_i)$ for $i\le m$. Then
$$
F_{i-1}\cup F_i \subset F_i^*.
$$
  Set
$s=(4c_2)^{-1}r^2$. Then $ s\asymp \rho(z_i, z_{i-1})^2\asymp
\rho_{z_i}( r_i)^2$. As (\ref{e:1236wq}), we have
\begin{equation}\label{e:71613} p_s(y',x')\ge \frac{c_4}{\nu(F^*_i)} \hbox{~~ for~ } y'\in
F_{i-1}, \, x'\in F_i. \end{equation} By Lemma \ref{l:1331}, for
$y'\in F_{i-1}$,
$$
\P_{y'}(X_s\in F_i)\ge c_4\frac{\nu(F_i)}{\nu(F_i^*)}\ge c_5.
$$
 Note that $$t-ms=t-m\cdot (4c_2)^{-1}r^2\ge t-(2c_2\rho(x,y)/r)\cdot (4c_2)^{-1}r\cdot(t/\rho(x,y))= \tfrac{t}{2}.$$
So, as (\ref{e:1236wq}) we can get
\begin{equation}\label{e:71614}
p_{t-ms}(x,y')\ge \frac{c_{6}}{V_{\rho}(x,t^{1/2})}~~{\rm~for~}y'\in
F_m.
\end{equation}
Therefore,
\begin{align*}
p_t(x,y)&=p_t(y,x)\ge \mu_x^{-1}\P_y(X_{is}\in F_{i}, 1\le i\le m,
X_{t}=x)\\
&\ge  c_5^{m}  \min_{y'\in F_{m}}\P_{y'}(X_{t-ms}=x)\mu_x^{-1}\\
&=c_5^{m}  \min_{y'\in F_{m}}p_{t-ms}(x,y')
\ge c_5^{m} \frac{c_6}{V_{\rho}(x,t^{1/2})}\\
 &\ge  \frac{c_6}{V_{\rho}(x,t^{1/2})}\exp\{-c_5'm\}
 \ge
\frac{c_6}{V_{\rho}(x,t^{1/2})}\exp\{-2c_2c_5'\frac{\rho(x,y)^2}{t}\},
\end{align*}
which implies (\ref{lower:1608w}). We have completed the proof of
Theorem \ref{T:lower1848}. \qed

\section{Proof of Theorem \ref{t:main}}
{\it Proof of Theorem \ref{t:main}.} (1) By Theorem
\ref{T:lower1848} and Lemma \ref{l:0900}, if $\al<d-2$ then
\begin{align*}
\int_1^\infty p_t(0,0)dt\le  c\int_1^\infty
t^{-(d+\al)/(2+2\al)}dt=\frac{2+2\al}{d-2-\al}c<\infty.
\end{align*}
Hence if $\al<d-2$ then $X$ is transient. Similarly, if  $\al\ge
d-2$ then $\int_1^\infty p_t(0,0)dt=\infty$ and so $X$ is
recurrent.\\

(2)  Let  $X'$ be an independent copy of $X$. We use $\P_{x,x'}$ for
the probability measure of the processes $X$ and
 $X'$ which start from $x$ and $x'$ respectively.

 If $d=1$ then
\begin{align*}
 \int_1^\infty \P_{0,0}(X_t=X_t'=0)dt=& \int_1^\infty
 \P_{0}(X_t=0)\P_{0}(X'_t=0)dt\\
 =&\int_1^\infty
 p_t(0,0)^2dt\ge c\int_1^\infty
t^{-2(1+\al)/(2+2\al)}dt=\infty.
 \end{align*}
So, $(X,X')$ is recurrent, which implies $X$ and $X'$ collide at the
origin infinitely often.

 Let  $d=2$. Fix $\lambda= \lceil 100c_{\ref{t:09051}.1}\rceil> 100$. For
$k\ge 1$, we set
$$
t_k={\lambda}^{2k(1+\al)},
$$
$$
\mathbb{T}_k=B(0, 2{\lambda}^{k})-B(0,{\lambda}^k),
$$
 $$ \theta_k=\inf\{t\ge 0:
|X_t|\ge{\lambda}^{k+1}\},~~~~\theta_k'=\inf\{t\ge 0: |X'_t|\ge
{\lambda}^{k+1}\}
$$
and
$$
H_k=\int_{0}^{\theta_{k}\wedge \theta_k^{'}\wedge
2t_k}1_{\{X_t=X_t'\in \mathbb{T}_k\}}dt.
$$
So, if $H_k>0$ then there exists at least one collision of  $X$ and
$X'$  before their breaking out of $B(0, \lambda^{k+1})$.  We shall
use the second moment method to estimate the probability of the
event $\{H_k>0\}$ as the approach of  \cite{CCD,CCD2}. Fix $x,y\in
B(0,{\lambda}^k)$. Then
\begin{align}
\nonumber\E_{x,y}(H_k)= &\int_{0}^{2t_k} \P_{x,y}(X_t=X'_t\in
\mathbb{T}_k,
\theta_k>t,\theta_k^{'}>t)dt\\
\nonumber\ge&\int_{t_k}^{2t_k} \sum_{u\in
\mathbb{T}_k}\P_{x,y}(X_t=u, X'_t=u,
\theta_k>t, \theta'_k>t)dt\\
\label{e:1733w}=&\int_{t_k}^{2t_k} \sum_{u\in
\mathbb{T}_k}\P_{x}(X_t=u, \theta_k>t)\P_y(X'_t=u,\theta'_k>t)dt.
\end{align}
Note that $t_k= {2}^{-2-2\al}\rho_0(2{\lambda}^{k})^2$ and
$\lambda^{k+1}=\lceil100 c_{\ref{t:09051}.1}\rceil\lambda^k $.
 Employing Theorem
\ref{t:09051} on $B(0, 2{\lambda}^{k})$, we get for each $u, v\in
B(0, 2{\lambda}^{k})$ and $t\in [t_k, 2t_k]$,
$$
\P_{u}(X_t=v, \theta_k>t)\ge \frac{c\nu_v}{V(0, 2\lambda^k)}.
$$
By  Lemma \ref{l:0900},   for    $v\in \mathbb{T}_k$ ,
$$
\frac{\nu_v}{V(0, 2\lambda^k)}\ge
c\frac{|v|^\al}{(2\lambda^k)^{2+\al}}\ge c'\lambda^{-2k}.
$$
Hence $\P_u(X_t=v,\theta_k>t)\ge c\lambda^{-2k}$ for each $u\in
\{x,y\}, v\in \mathbb{T}_k$ and $t\in [t_k, 2t_k]$. Therefore,
inequality  (\ref{e:1733w}) becomes
\begin{equation}\label{e:1738w}
\E_{x,y}(H_k)\ge  (c\lambda^{-2k})^2\cdot|\mathbb{T}_k|\cdot t_k\ge
c^2 \lambda^{-4k}\cdot c'(\lambda^k)^2\cdot \lambda^{2k(1+\al)}=c''
{\lambda}^{2k\al}.
\end{equation}

On the other hand, for any $u\in \mathbb{T}_k$,
\begin{align}
\nonumber\E_{u,u}(H_k)\le &\int_{0}^{2t_k}\sum_{w\in \mathbb{T}_k} [\P_{u}(X_t=w)]^2dt\\
\nonumber\le &\frac{\max_{w\in
\mathbb{T}_k}\nu_w}{\nu_u}\int_{0}^{2t_k} \sum_{w\in
\mathbb{T}_k}\P_{u}(X_t=w)\P_w(X_t=u)dt\\
\nonumber\le& c\int_{0}^{2t_k} \P_{u}(X_{2t}=u)dt\le
c\nu_u^2+c\int_{\nu_u^2}^{2t_k}
\P_{u}(X_{2t}=u)dt\\
\nonumber\le& c\nu_u^2+\int_{\nu_u^2}^{2t_k}
\frac{c'\nu_u}{V_{\rho}(u, t^{1/2})}dt  \\
\nonumber\le&  c\nu_u^2+c''\nu_u^2\int_{\nu_u^2}^{2t_k} t^{-1}dt,
\end{align}
where the last second inequality is by (\ref{e:1845w}), while the
last by Lemma \ref{l:0900}. Hence
\begin{equation}\label{e:1840w}
\E_{u,u}(H_k)\le c\nu_u^2(1+\log (2t_k)-\log(\nu_u^2))\le c'
\lambda^{2k\al}\cdot (\log
(\lambda^{2k(1+\al)})-\log(\lambda^{2k\al}))= c''k\lambda^{2k\al}.
\end{equation}
By the strong Markov property,
\begin{align*}
\E_{x,y}(H_k^2)=&2\E_{x,y}\left(\int_{0}^{\theta_{k}\wedge
\theta_k^{'}\wedge 2t_k}1_{\{X_t=X'_t\in
\mathbb{T}_k\}}dt\int_{t}^{\theta_{k}\wedge \theta_k^{'}\wedge
2t_k}1_{\{X_s=X'_s\in \mathbb{T}_k\}}ds\right)\\
\le &2\E_{x,y}\left(\int_{0}^{\theta_{k}\wedge \theta_k^{'}\wedge
2t_k}1_{\{X_t=X'_t\in \mathbb{T}_k\}}\E_{X_t,X'_t}(H_k)dt\right)\\
\le&2\sup_{u\in \mathbb{T}_k} \E_{u,u}(H_k) \E_{x,y}(H_k).
\end{align*}
So, by (\ref{e:1840w}), (\ref{e:1738w}) and the Cauchy Schwarz
inequality,
$$
\P_{x,y}(H_k>0)\ge \frac{[\E_{x,y}(H_k)]^2}{\E_{x,y}(H_k^2)}\ge
\frac{\E_{x,y}(H_k)}{2\sup_{u\in \mathbb{T}_k} \E_{u,u}(H_k)}\ge
\frac{c}{k}.
$$
Therefore, when $X$ and $X'$ start from $x,y\in B(0,\lambda^k)$
respectively, the probability that they will collide  before their
breaking out $B(0,\lambda^{k+1})$, is not less than $\tfrac{c}{k}$.
 Note that $\sum_k\tfrac1k=\infty$. Using the second
Borel-Cantelli Lemma as \cite[Theorem 1.1]{CCD2}, we  prove that $X$
and $X'$ collide infinitely often when $d=2$.
\\


(3) Let $d\ge 3$. For $k\ge 0$, set
$$
\mathbb{T}_k=B(0,2^{k+1})-B(0,2^k)~~~{\rm
and~~~}Z_{k}=\int_{0}^\infty 1_{\{X_t=X_t'\in \mathbb{T}_k\}}dt.
$$
Then
\begin{align*}
\E_{0,0}(Z_k)
=&\sum_{u\in \mathbb{T}_k}\int_{0}^\infty [\P_0(X_t=u)]^2dt\\
=&\sum_{u\in \mathbb{T}_k}\int_{t_k}^\infty [\P_0(X_t=u)]^2dt+\sum_{u\in \mathbb{T}_k}\int_{s_k}^{t_k} [\P_0(X_t=u)]^2dt+\sum_{u\in \mathbb{T}_k}\int_0^{s_k} [\P_0(X_t=u)]^2dt\\
=&I_1+I_2+ I_3,
\end{align*}
where $ s_k=(1\vee 2^{k\al})2^{k(1+\al)}$ and $t_k=2^{k(2+2\al)}.$
We shall deal with the three sums separately. Since $t_k\ge
c\rho(0,u)^2$ for $u\in \bT_k$, we can use Theorem \ref{T:lower1848}
and  Lemma \ref{l:0900}, and get
\begin{align*}
 I_1\le& \sum_{u\in
\mathbb{T}_k}\int_{t_k}^\infty
\frac{c\nu_u^2}{V_{\rho}(0,t^{1/2})V_{\rho}(u, t^{1/2})} dt\\
 \le& |\mathbb{T}_k|\cdot\max_{u\in \mathbb{T}_k}\nu_u^2\cdot  \int_{t_k}^\infty
c't^{-(d+\al)/(1+\al)}dt\\
\le& 2^{dk}\cdot c'' 2^{2k\al}\cdot c'''
(2^{2k(1+\al)})^{1-(d+\al)/(1+\al)}\\
=&c2^{k(2+2\al-d)}.
\end{align*}
Next, since $(1\vee \nu_u)\rho(0,u)\ge cs_k$ and $t_k^{1/2}\le
c'|u|^{1+\al}$ for $u\in \bT_k$, using Theorem \ref{T:lower1848} and
Lemma \ref{l:0900} again gives
\begin{align*}
 I_2\le& \sum_{u\in
\mathbb{T}_k}\int_{s_k}^{t_k}
\frac{c\nu_u^2}{V_{\rho}(0,t^{1/2})V_{\rho}(u, t^{1/2})}\exp\left(-\frac{\rho(0,u)^2}{ct}   \right) dt\\
\le&\sum_{u\in \mathbb{T}_k}\int_{s_k}^{t_k}
\frac{c'\nu_u^2}{t^{(d+\al)/(2+2\al)}\cdot t^{d/2}|u|^{-(d-1)\al}}\exp\left(-\frac{2^{2k(1+\al)}}{c't}\right) dt\\
\le& |\mathbb{T}_k|\cdot\max_{u\in
\mathbb{T}_k}\{\nu_u^{2}|u|^{(d-1)\al}\}\cdot \int_{0}^{\infty}
c't^{-(d+\al)/(2+2\al)-d/2}\exp\left(-\frac{2^{2k(1+\al)}}{c't}\right)dt\\
\le&2^{dk} \cdot c''2^{(d+1)k\al}\cdot
(2^{2k(1+\al)})^{1-(d+\al)/(2+2\al)-d/2}\cdot c'''\int_0^\infty x^{(d+\al)/(2+2\al)+d/2-2}e^{-x}dx\\
=&2^{k(2+2\al-d)}\cdot c''c'''\int_0^\infty
x^{(d+\al)/(2+2\al)+d/2-2}e^{-x}dx.
\end{align*}
Since $d\ge 3$, we have $\int_0^\infty
x^{(d+\al)/(2+2\al)+d/2-2}e^{-x}dx<\infty$ and so,
$$
I_2\le c 2^{k(2+2\al-d)}.
$$
For the remaining term, applying  Theorem \ref{T:lower1848} we still
have
\begin{align*}
 I_3\le& \sum_{u\in
\mathbb{T}_k}\int_{0}^{s_k}
({\nu}_{u}/{\nu}_0)\exp\left(-c(\nu_0\vee \nu_u)^{-1}
\rho(0,u)\Big(1\vee\log\big( (\nu_0\vee\nu_u) \rho(0,u)/t\big)\Big)\right) dt\\
\le& |\mathbb{T}_k|\cdot s_k \cdot\max_{u\in \mathbb{T}_k}\nu_u
\exp\big(-c(\nu_0\vee \nu_u)^{-1}
\rho(0,u)\big)\\
\le& 2^{dk}\cdot (1\vee 2^{k\al})2^{k(1+\al)} \cdot 2^{k\al} \exp\left(-c' (1\vee 2^{k\al})^{-1}\cdot 2^{k(1+\al)}\right)\\
=&2^{ck}e^{-c'2^{k(1+\al\wedge 0)}} \le c''2^{k(2+2\al-d)}.
\end{align*}
Therefore,
\begin{equation}\label{e:0730w}
\E_{0,0}(Z_k)\le c2^{k(2+2\al-d)}.
\end{equation}

On the other hand, once $X$ and $X'$ collide at some vertex $u$ and
some time $t$,  then with at least $e^{-2}$ probability   they will
stick together  during   time $[t,t+\nu_u/\mu_u)$, which implies
$$\E_{0,0}(Z_k|Z_k>0)\ge c\min_{u\in \mathbb{T}_k}\nu_u/\mu_u\ge
c'2^{2k\al}.$$ So, for each $k\ge 0$,
$$
\P_{0,0}(Z_k>0)= \frac{\E_{0,0}(Z_k)}{\E_{0,0}(Z_k|Z_k>0)}\le
c2^{-(d-2)k}.
$$
Therefore,
$$
\sum_k\P_{0,0}(Z_k>0)\le c\sum_k 2^{-(d-2)k}<\infty.
$$
By the Borel-Cantelli Lemma, we completed the proof of (3).\qed

\noindent {\bf Xinxing Chen}

\noindent Department of Mathematics, Shanghai Jiaotong University,
Shanghai, China, 200240.

\noindent E-mail: {\tt chenxinx@sjtu.edu.cn}

\end{document}